\documentclass[10pt,twocolumn]{article}
\usepackage[affil-it]{authblk}
\usepackage{natbib}
\setcitestyle{round}
\usepackage[font=scriptsize]{caption}
\usepackage[dvips]{graphicx}
\graphicspath{{./}}
\DeclareGraphicsExtensions{.eps}

\usepackage[cmex10]{amsmath}
   \title{An adjustable-width window with good dynamic range.}

\author[1]{\rm I. M. Stewart}
\affil[1]{Argelander Institut f\"ur Astronomie, Universit\"at Bonn,
             Auf dem H\"ugel 71, 53121 Bonn, Germany}
\begin{document}
   \maketitle

   \abstract{
A new variable-width window is presented and compared with several other windows, both of variable and fixed widths. The comparison focuses on sensitivity and dynamic range. The equivalent noise bandwidth or ENBW (or rather, its reciprocal) is used as a proxy for the first; maximum sidelobe level and high-frequency roll-off in the Fourier transform, for the second. The new window can access any value of ENBW by appropriate choice of the width parameter. At any given value of ENBW below about 3, a setting can be found at which the sidelobes of the window are lower than those of any other in the moderate frequency regime below about 100 cycles. 
   }

%++++++++++++++++++++++++++++++++++++++++++++++++++++++++++++
\section{Introduction}

A commonly encountered class of signals consists of one or more pure sinusoidal tones plus noise. For such signals, the Fourier transform (FT) is often the analysis method of first choice because it concentrates the power of each tone into a delta function located at the tone frequency, and also displays the correlation properties of the noise. Unfortunately it is not possible to compute integrals of continuous functions which extend to infinity. In actual computations therefore one must calculate an approximation to the Fourier transform which will necessarily involve a finite number $N$ of samples of the function. The input signal is notionally multiplied by a uniform window extending over the interval $[-T/2,T/2]$ (which is sometimes called a top hat window because of its shape), then sampled within that interval at constant intervals $\Delta t = T/N$. The delta function for each tone becomes convolved with the Fourier transform of the top hat window, which is a sinc function.

If the Fourier transform of the windowed and sampled input is also sampled at $N$ intervals separated by $\Delta f = 1/T$, the result is known as the Discrete Fourier Transform or DFT \citep{bracewell_2000}. If the sinusoids contained in the original signal happen to have periods $T_j$ which are integer fractions of $T$, then the DFT will concentrate all the power of such tones into a single bin of the output. In effect the samples occur at the zeros of the sinc function. In the general case, however, non-periodicity of the sinusoids at the window edges is to be expected: this will generate Fourier components which were not there in the original signal. The effect is to redistribute power out of the bin corresponding to the true frequency of the tone into adjacent bins. This `leakage' can make it difficult to find other tones close to and/or fainter than the strongest.

A classic way to ameliorate this problem (nothing will make it entirely disappear) is to choose a different window function. A window which tapers gradually toward zero as $|t| \to T/2$ is usually found to cause much less leakage than one which drops abruptly to zero at the boundaries. A variety of functions having differing properties have been proposed over the years. Commonly, although not always, a window function will be real-valued, symmetrical about $t=0$, and everywhere nonnegative, the value decreasing monotonically from a maximum at $t=0$ to a minimum at $|t|=T/2$. Only windows having these properties are discussed in the present paper. Such a window has a Fourier transform which typically is reminiscent of a sinc function, namely with a relatively broad central positive-valued lobe surrounded by oscillations, known as sidelobes, which decrease steadily in amplitude with increasing frequency $f$. Often too (as with the sinc function) the rate of decrease of this amplitude is proportional to $f$ raised to some negative power.

\citet{harris_1978} defined several criteria for assessing the performance of windows and tabulated the respective values for many examples. In the present paper three only of these criteria will be considered: the equivalent noise bandwidth or ENBW, the amplitude of the largest sidelobe, and the far-field rate of decrease of sidelobe amplitude (sometimes known as the roll-off). These quantities are described in more detail in section \ref{ss_criteria}.

The present paper describes a new window function which includes a parameter which allows one to vary the ENBW over its entire possible range. For convenience, and for a reason which is made clear in section \ref{ss_hyperbolic}, I refer to it as the hyperbolic window. The new window is compared here with a selection of others which excel in the properties of present interest, as described above. The hyperbolic window is not the only one of variable width, and the comparison makes it plain that one can often find an alternative which has a better value of either roll-off or sidelobe amplitude. However, as is shown in the present paper, one can still achieve roll-off rates of up to 60 decibels per octave with the hyperbolic window over a wide range of values of ENBW; and superior sidelobe levels can be attained in the near-field regime $Tf < \sim 100$, particularly in the range of ENBW values below 1.7, for which few well-performing windows are available. For this reason, and because of its flexibility, the hyperbolic window is suggested for applications in which sensitivity and near-field dynamic range need to be jointly optimized.

%++++++++++++++++++++++++++++++++++++++++++++++++++++++++++++
%\section{Windows and assessment criteria} \label{s_windows}
\section{Theory}
%oooooooooooooooooooooooooooooooooooooooooooooooooooooooooooo
\subsection{Generic windows} \label{ss_windows}

A general symmetric window function $w(|t|)$ is uniquely defined by its values in the interval $[0, T/2]$. Its Fourier transform $W(f)$ is thus clearly
\begin{displaymath}
	W(f) = 2 \int_0^{T/2} dt \ w(t) \; cos(2 \pi f t).
\end{displaymath}
Throughout the present paper I assume that $w$ has been normalized such that $W(0)=1$.

In order to do $N$-point DFT calculations with a window, one has to sample $w(t)$ at $1+N/2$ locations, the samples being defined as
\begin{displaymath}
	w_j = w(j T / N) \textrm{ for } j \in \{0, \ \dots, \ N/2\}.
\end{displaymath}
(I assume throughout that $N$ is even.) Conventionally then one stores the symmetric samples corresponding to $t < 0$ in the upper part of the vector according to the prescription
\begin{displaymath}
	w_j = w_{N-j} \textrm{ for } j \in \{N/2+1, \ \dots, \ N-1\}.
\end{displaymath}
Note that $w_0$ and $w_{N/2}$ are centres of symmetry and thus each occurs only once in the sequence.

%**** mention here that w_j can be decomposed into a sum of cosines?

The Fourier transform of the continuous function $w$ can be more closely approximated via the DFT, i.e. with a finer frequency spacing, if the technique of zero-padding is used \citep{bracewell_2000}. Suppose one wished to sample the FT with spacing $\delta f = 1/a$. (Note that $a$ must be $>1$ and that $M = aN$ must be an integer.) The procedure is to construct a new vector $w^\prime_j$ of length $M = aN$ for which
\begin{displaymath}
	w^\prime_j = w_j \textrm{ for } j \in \{0, \ \dots, \ N/2\}
\end{displaymath}
and
\begin{displaymath}
	w^\prime_j = w_{j-M+N} \textrm{ for } j \in \{M+1-N/2, \ \dots, \ M-1\},
\end{displaymath}
the remainder being filled with zeros.

%oooooooooooooooooooooooooooooooooooooooooooooooooooooooooooo
\subsection{Assessment criteria} \label{ss_criteria}

The three criteria listed by \citet{harris_1978} used in the present paper are the equivalent noise bandwidth (ENBW), the maximum sidelobe power, and the roll-off rate. The ENBW is defined as
%\begin{displaymath}
\begin{equation} \label{equ_enbw}
	ENBW = N \frac{\sum_{j=0}^{N-1} w^2_j}{\left( \sum_{j=0}^{N-1} w_j \right)^2}.
\end{equation}
%\end{displaymath}
This gives the ratio between the gain for incoherent sources versus coherent ones, and thus is the reciprocal of sensitivity, which can be described as the detectability of a pure tone (i.e. coherent) among noise (incoherent). The smallest possible value of ENBW (corresponding to the best sensitivity) is unity, which can be attained only by the top-hat window. As can be seen from table 1 of \citeauthor{harris_1978}, ENBW correlates well with measures of the width of the central lobe of the window FT $W(f)$, but the exact proportionality varies from window to window.

The `sidelobe height' is defined here as the maximum value of $W^2(f)$ for $f > f_0$, where $f_0$ is the frequency of the first zero-crossing of $W$. Note that the FTs of some windows don't cross zero at all, or at least not until very large $f$; thus a sidelobe height cannot meaningfully be defined for them. In the present paper, following conventional practice, sidelobe heights are presented in decibels.

The FT of any function $w$ which is only non-zero within a finite interval, and analytic within that interval, is found to decrease in amplitude asymptotically at some finite order $(n+1)/2$ of frequency, the value of $n$ being determined by the order of the first differential of $w$ which is non-zero within or at the interval boundaries. On a logarithmic scale, a change in FT envelope amplitude proportional to frequency raised to the power $-(n+1)/2$ is equivalent to about $6 \times (n+1)$ decibels decrease per octave. I find it simpler just to quote $-(n+1)$: hence all roll-off values quoted in the present paper follow this scheme. The top hat window for example has, in these terms, a roll-off of $-1$.

In computational practice the maximum frequency obtained is limited to $N/(2T)$, $N$ being as before the number of samples of $w$. This may not be large enough to allow $W$ to reach its asymptotic value. It should also be pointed out that the finite precision of computation places an ultimate limit on dynamic range, so for purposes of optimizing this over the widest range of accessible frequencies, one needs to look farther than the theoretical asymptotic roll-off, and also ask how soon the asymptotic behaviour becomes manifest.

Both sidelobe height and roll-off may be used as proxies for dynamic range, since they both bear on the detectability of secondary tones in the presence of the strongest.

%oooooooooooooooooooooooooooooooooooooooooooooooooooooooooooo
\subsection{The hyperbolic window} \label{ss_hyperbolic}

The hyperbolic window is defined for $|t|<T/2$ as
\begin{displaymath}
%%%	w(t) = 0.5 \{1 + \cos[\pi z(s,t)]\}
%%	w(t) = \cos^{2 \alpha}[z(s,t)]
%	w(t) = \cos^{2 \alpha}[\pi z(s,t) / 2]
	w(t) = \cos^{2 \alpha} \left[ \frac{\pi z(s,t)}{2} \right]
\end{displaymath}
where
\begin{displaymath}
%	z(s,t) = \frac{2 |t| (1 - s) / T}{1 - s (4 |t| / T - 1)}.
	z(s,t) = \frac{\tau (1 - s)}{1 - s (2 \tau - 1)}
\end{displaymath}
for
\begin{displaymath}
	\tau = 2 |t| / T.
\end{displaymath}
The variable width is achieved via the warp factor $s$, for which the useful range of variation is -1 to 1. As $s \to -1$, $w$ becomes infinitely narrow; at $s=0$, $z=\tau$, thus $w$ becomes identical to a simple cosine window, raised to the power $2 \alpha$; as $s \to 1$, $w$ approaches the top hat. Note e.g. that for $s=0$ and $\alpha=1$, the hyperbolic and Hann windows \citep{harris_1978} are identical.

No closed-form expression is known for the general Fourier transform.

The name `hyperbolic' was chosen because $z$ can be rearranged to the form
\begin{displaymath}
	z(s,t) = A(s) + \frac{B(s)}{C(s) + |t|}.
\end{displaymath}

The Fourier transform of the hyperbolic window undergoes a marked changes as the warp factor $s$ departs from the `Hann' value of zero. For infinitesimal values of $s$ it is easy to show via Taylor expansion that the difference $\Delta w$ between the hyperbolic and Hann windows is (for integer $\alpha$) equal to
%\begin{equation} \label{equ_perturbation}
%%%	w_\mathrm{hyperbolic} - w_\mathrm{Hann} = 2 \pi \alpha s \cos^{2 \alpha - 1} \left( \pi \frac{|t|}{T} \right) \sin \left( \pi \frac{|t|}{T} \right) \frac{2|t|}{T} \left( 1 - \frac{2|t|}{T} \right).
%%	\Delta w = 2 \pi \alpha s \cos^{2 \alpha - 1} \left( \pi \frac{|t|}{T} \right) \sin \left( \pi \frac{|t|}{T} \right) \frac{2|t|}{T} \left( 1 - \frac{2|t|}{T} \right).
%	\Delta w = \frac{\pi \alpha s}{2^{2\alpha-2}} \tau ( 1 - \tau ) \sum_{k=0}^{\alpha-1} \binom{2\alpha}{k} \sin [(\alpha - k) \pi \tau ]
%\end{equation}
\setlength{\arraycolsep}{0.0em}
\begin{eqnarray} \label{equ_perturbation}
	\lefteqn{\Delta w = 2^{(2-2\alpha)} \pi \alpha s \tau ( 1 - \tau ) \times {}}\nonumber\\
	& & {\times}\: \sum_{k=0}^{\alpha-1} \binom{2\alpha}{k} \frac{\alpha-k}{\alpha} \sin [(\alpha - k) \pi \tau ]
%	& & {} \times \sum_{k=0}^{\alpha-1} \binom{2\alpha}{k} \frac{\alpha-k}{\alpha} \sin [(\alpha - k) \pi \tau ]
\end{eqnarray}
\setlength{\arraycolsep}{5pt}%
for $0 \le \tau \le 1$. Since this begins to have discontinuities at derivative order 3, one expects it to have a roll-off of -4. In fact this is what is observed (See for example Fig. \ref{fig_f}). However, the Fourier transform of equation \ref{equ_perturbation} appears to consist of a smooth or DC component of this roll-off value plus an oscillatory component which, particularly at high values of $\alpha$, exhibits a much steeper roll-off. The upshot is that the total FT for the hyperbolic window tends to resemble that of the Hann window raised to the same $\alpha$ exponent, plus a DC component at the much shallower roll-off of $-4$. Clearly the shallow-slope DC component is going to dominate the total at high frequencies; at low, it may or may not remain the dominant component, depending on the values of $s$ and $\alpha$.

%++++++++++++++++++++++++++++++++++++++++++++++++++++++++++++
\section{Comparisons} \label{s_comparisons}
%oooooooooooooooooooooooooooooooooooooooooooooooooooooooooooo
\subsection{The contenders} \label{ss_contenders}

The highest-sidelobe and roll-off for the hyperbolic window for several samples of values of $s$ and $\alpha$ are compared here against matching values for a selection of other window types, both variable and non-variable. These are described in the following lists.

Two variable-width windows were included:
\begin{itemize}
  \item The Tukey window \citep{harris_1978, tukey_1967} is a simple modification of the Hann window to allow ENBW to be varied from the Hann value to unity while retaining the $-3$ roll-off of the Hann. One could also choose to raise the Tukey to some arbitrary power, but in fact this does not improve the roll-off, while at the same time worsening the sidelobes. The lack of improvement results because the 2nd-order and higher derivatives remain discontinuous at the boundary between the constant part and the cosine part of this window. The FT of the Tukey looks rather uneven and exhibits large beats as the window width approaches either end of its range, due to the fact that it comprises two different functions spliced together.

  \item The Planck window \citep{mckechan_2010} has continuity of derivatives of all orders, everywhere. Because of that it has an asymptotic roll-off which increases without bound. Beats also occur in the FT of this window, for the same reason as for the Tukey window.
\end{itemize}

The fixed-width windows are:
\begin{itemize}
  \item The top hat window is included because it forms the end point for the transformation of all three variable-width windows.

  \item Hann \citep{harris_1978}: one of the simplest, yet it has reasonable values of all three criteria considered here. It is also a shape reachable by both the hyperbolic and Tukey windows. Squared and cubed Hann windows are also compared here.

  \item De la Vall\'{e}-Poussin \citep{parzen_1961}: included again because of its good compromise of values.

  \item Bohman \citep{bohman_1960}: properties similar to the De la Vall\'{e}-Poussin, and to the squared Hann.

  \item Kaiser \citep{kaiser_1980}: included because it has excellent sidelobe suppression at moderate width, although it rolls off only as $1/f$. It is a variable-width formulation but, following \citet{nuttall_1981}, only the window normalized by $I_0(3\pi)$ is considered here.

  \item Nuttall \citep{nuttall_1981}: Nuttall explored several windows consisting of sums of 3 or 4 cosine terms, optimized mostly for sidelobe level, but also at respectable roll-off values. Three are included for comparison here: namely the 3-term window corresponding to his Fig. 9, which is labeled here `Nuttall 3'; and 4-term windows `Nuttall 4a' and `Nuttall 4b', corresponding to his Figs. 11 and 12 respectively.
\end{itemize}

Many more fixed-width windows have been described, but most of them have poor roll-off, or too-large ENBW, so do not commend themselves for high dynamic range, high sensitivity work. In any case I do not pretend to make here a comprehensive survey of windows.

Note also that, except for the Hann, I have not thoroughly explored the changes in the fixed-width windows obtainable by raising them to higher powers. The few trials which have been made hold out little promise. The Nuttall windows for example are carefully balanced to yield optimum sidelobe levels; raising them to a higher exponent destroys this balance, with resulting rapid deterioration of the sidelobe performance. The same is observed when the hyperbolic transform is tried with these windows.

%oooooooooooooooooooooooooooooooooooooooooooooooooooooooooooo
\subsection{Results} \label{ss_graphs}

For all the results quoted in the present section, window functions were sampled at $N=4096$ or $2^{12}$ points. Since a zero pad multiple of 16 was used, the discrete Fourier transform operated on $2^{16}$ points in total.

The ENBW is easily calculated from equation \ref{equ_enbw}. The highest sidelobe level is best estimated via a DFT of a sufficiently zero-padded version of the sampled window. Roll-off slopes were found the most difficult to estimate in practice. For the present paper, these were estimated from log-log scaled, squared FT curves in an iterative process in which an initial guess at the slope of the FT envelope was progressively refined. This slope was estimated over a range of values of frequency around $150/T$.

For the machine used to perform the calculations, floating-point precision was reported as 2.2e-16. This is consistent with the observed noise floor in the Fourier transforms, which lies at about $-340$ decibels below the peak power. Only the Planck and hyperbolic FTs routinely reach this noise floor before $Tf = N/2$.

\begin{figure}[!t]
	\centering
	\includegraphics[width=2.5in]{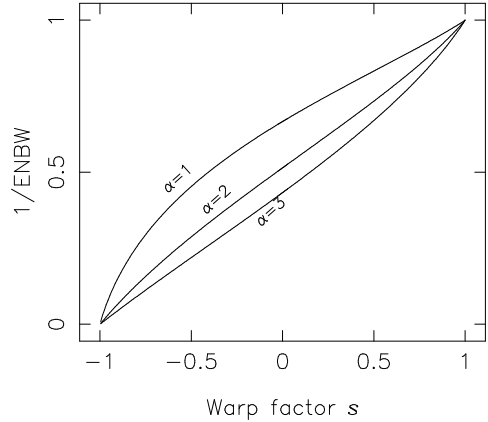}
	\caption{Variation of the equivalent noise bandwidth (ENBW) with `warp factor' $s$ for the hyperbolic window. Three curves are provided for different values of the exponent $\alpha$. Note that the reciprocal ENBW is in fact what is plotted, so as to display its entire range.}
	\label{fig_g}
	\end{figure}

The warp parameter $s$ of the hyperbolic window is monotonically related to the reciprocal of ENBW. Fig. \ref{fig_g} displays this relation for three values of $\alpha$.

\begin{figure}[!t]
	\centering
	\includegraphics[width=2.5in]{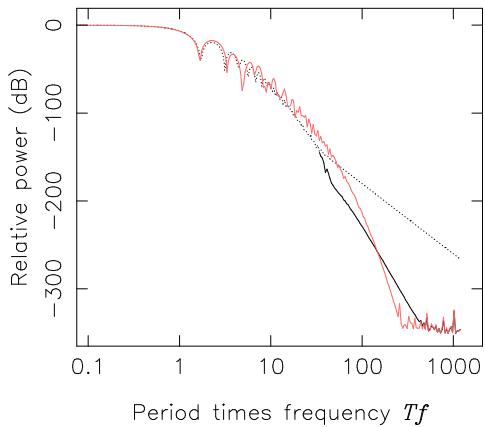}
	\caption{The Fourier power spectra of the hyperbolic (dotted and solid black lines) and Planck (solid red line) windows are shown. The value of ENBW is 1.46 for both. The Planck window was calculated with an $\epsilon$ value of 0.4; for the hyperbolic window, $\alpha=4$ and $s=0.606$. The dotted curve plots the full FT of the hyperbolic window, whereas the solid black line traces the oscillatory residual. Note also the noise floor at about $-340$ decibels due to the finite computational precision.}
	\label{fig_j}
	\end{figure}

In section \ref{ss_hyperbolic} it is described how the FT of the hyperbolic window often exhibits a DC component of roll-off about $-4$ plus a significantly steeper oscillatory component. An example of the full hyperbolic window FT is shown in Fig. \ref{fig_j}, compared to a Planck window of the same value of ENBW. The separation between the DC and oscillatory components is clearly visible at frequencies greater than about $Tf \sim 30$. At frequencies below this value, the total FT is oscillatory - it passes through zero - even though this is not readily apparent from Fig. \ref{fig_j}, because of the plotting algorithm.

The tendency of the frequency at which the DC component ceases to dominate the whole is to move toward lower values as the ENBW increases. Past a certain value of ENBW, which depends on the value of the exponent $\alpha$, the DC component is dominant over the entire spectrum; in other words, the FT ceases to be oscillatory at any frequency. Shortly beyond this critical value of ENBW, local maxima are no longer found in the near-field transform. These critical points can be observed in Fig. \ref{fig_d} as the places beyond which sidelobe height as defined in section \ref{ss_criteria} ceases to be meaningful.
%An example comparison between hyperbolic and Planck windows of the same ENBW value is shown in Fig. \ref{fig_j}. Here too the `residual' oscillatory part of the hyperbolic window FT has been plotted as well as its total FT. As an example of the dynamic range obtainable at mid-frequency values by the respective windows, see Fig. \ref{fig_k}. Here a secondary tone at frequency $=100/T$ and of amplitude $10^{-10}$ smaller than the primary has been added before Fourier transforming. The secondary peak can be seen in both window FTs, but whereas one would probably not care to go much fainter using the Planck window, there is ample room to detect a 5-fold fainter tone in the hyperbolic residuals.

\begin{figure}[!t]
	\centering
	\includegraphics[width=2.5in]{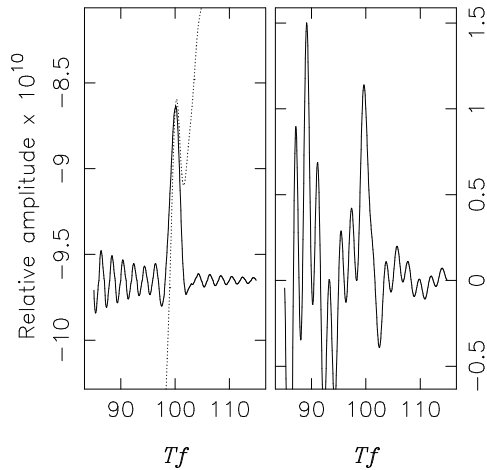}
	\caption{The same two window FTs as in Fig. \ref{fig_j} are plotted here again over a smaller range of frequencies, and with a secondary pure tone added at $Tf=100$ and of amplitude $10^{-10}$ relative to the $f=0$ value. The hyperbolic window is shown on the left, the Planck window on the right. In the LH panel, the dotted line shows the full FT, the solid line shows the oscillatory residual. The vertical range is the same in both panels but the full-FT plot for the hyperbolic has been offset vertically.}
	\label{fig_k}
	\end{figure}

In practice it is easy to estimate and remove the DC component via interpolation of a non-padded DFT of the window. Effectively then, as far as detecting weak secondary tones is concerned, the limiting factor is the oscillatory `residual' component of the FT, rather than its total amplitude. This is demonstrated in Fig. \ref{fig_k}, which shows a zoomed-in portion of the spectra in Fig. \ref{fig_j}. Here also a secondary tone at frequency $=100/T$ and of amplitude $10^{-10}$ smaller than the primary has been added before Fourier transforming. The secondary peak can be seen in both window FTs, but whereas one would probably not care to go much fainter using the Planck window, there is ample room at this frequency to detect a 5-fold fainter tone in the hyperbolic residuals.

At parameter settings for which the onset of DC dominance falls below the roll-off estimation frequency of about $150/T$, the roll-off of the hyperbolic window was estimated from the residual part of the FT; otherwise from the total FT.

\begin{figure}[!t]
	\centering
	\includegraphics[width=2.5in]{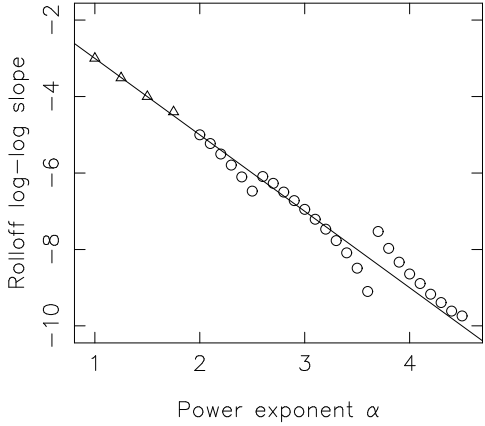}
	\caption{Roll-off of the hyperbolic window as a function of the exponent $\alpha$. For each value of $\alpha$, the warp factor $s$ is varied so as to maintain the ENBW at the Hann-window value of 1.5. The roll-off was estimated at frequency $f$ of $150/T$. An open triangle means the roll-off was estimated from the total curve; an open circle indicates the roll-off refers to the oscillatory residual. The difference between these is described in section \ref{ss_hyperbolic}. The solid line shows the value for the Hann window raised to the same exponent.}
	\label{fig_i}
	\end{figure}

\begin{figure}[!t]
	\centering
	\includegraphics[width=2.5in]{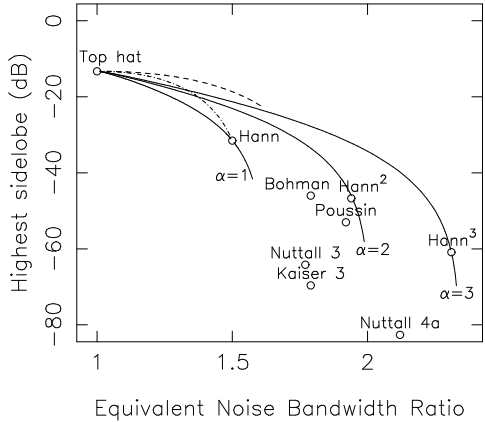}
	\caption{Plotted here for various windows is the ratio in decibels between the power of the highest sidelobe and the maximum power anywhere (which for all the present windows occurs at a frequency of zero). Circles indicate fixed-width windows. The dot-dashed line shows the Tukey window over its full permitted range of ENBW, the dashed likewise for the Planck window. The solid lines show the response of the hyperbolic window at three values of the exponent $\alpha$. Note that none of the hyperbolic-window curves extends over the full range of ENBW values accessible to the window. As explained in section \ref{ss_graphs}, this is because the Fourier transform of this window ceases to exhibit local maxima at a certain maximum value of ENBW (e.g. at ENBW $\sim 1.58$ for $\alpha=1$).}
	\label{fig_d}
	\end{figure}

Fig. \ref{fig_i} shows the estimated roll-off for the hyperbolic window at a number of values of exponent $\alpha$, with ENBW held at the Hann value of 1.5. It can be seen that the decrease is broadly in line with that obtained from raising the Hann window to the same power. Note though that raising the Hann window to exponents greater than unity also increases its ENBW beyond 1.5. The increase in roll-off steepness for the hyperbolic window can however be obtained \emph{without} sacrifice in ENBW, and thus sensitivity.

Fig. \ref{fig_d} compares the maximum sidelobe height for all the windows except the Nuttall 4b, which at $-93.3$ decibels lies off the bottom of the graph. One would assess the hyperbolic window as a relatively poor performer, except that its variable character allows exploration of areas of the parameter space not covered by any of the windows which offer better sidelobe levels at a small number of fixed (and relatively large) values of ENBW. Note also that compared to the other variable-width windows, the $\alpha=1$ hyperbolic is always better than the Tukey, at the same value ($-3$) of roll-off; and the $\alpha=3$ hyperbolic is always better than the Planck, at values of roll-off which, if not superior to the Planck, are certainly respectable (being steeper than for any of the fixed windows), and which, as is argued below, offer about the same practical performance.

%	\begin{figure}[!t]
%	\centering
%	\includegraphics[width=2.5in]{fig_d}
%	\caption{Plotted here for various windows is the ratio in decibels between the power of the highest sidelobe and the maximum power anywhere (which for all the present windows occurs at a frequency of zero). Circles indicate fixed-width windows. The dotted line shows the Tukey window over its full permitted range of ENBW, the dashed likewise for the Planck window. The solid lines show the response of the hyperbolic window at three values of the exponent $\alpha$. Note that none of the hyperbolic-window curves extends over the full range of ENBW values accessible to the window. As explained in section \ref{ss_graphs}, this is because the Fourier transform of this window ceases to exhibit local maxima at a certain maximum value of ENBW (e.g. at ENBW $\sim 1.58$ for $\alpha=1$).}
%	\label{fig_d}
%	\end{figure}

\begin{figure}[!t]
	\centering
	\includegraphics[width=2.5in]{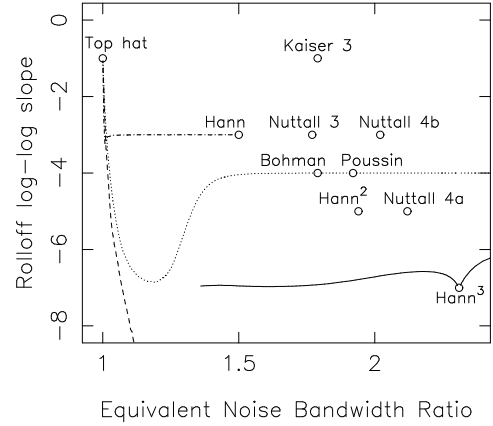}
	\caption{Plotted here is the roll-off slope of the hyperbolic window at $\alpha=3$, contrasted with that of the fixed-width windows (open circles), the Tukey (dot-dashed line) and the Planck (dashed line). The dotted line gives the roll-off calculated from the total curve; the solid line shows that calculated from the oscillatory residual. The difference between these is described in section \ref{ss_hyperbolic}. Present, but hard to see, is a discontinuity in the `total curve' line at the Hann-cubed value. As explained in section \ref{ss_hyperbolic}, this is due to the rapid growth of a shallow roll-off perturbation to the Hann-like FT as the warp factor moves away from zero.}
	\label{fig_f}
	\end{figure}

Fig. \ref{fig_f} compares roll-off for the $\alpha=3$ hyperbolic window against the other contenders. Although superior to most of the other windows, including the Tukey, which is stuck at a roll-off of $-3$ even when raised to a power, the best performer is clearly the Planck window, for which the present methods estimate a steepest roll-off value of about $-13$ at an ENBW just short of of 1.4. This number is actually not very meaningful, since the continuity of the Planck window at all orders of derivative means that its asymptotic roll-off is infinite. Indeed in plots of its FT (see Fig. \ref{fig_j} for example) the log-log slope is seen to increase monotonically with increasing $f$, unlike all other windows examined in the present study.

%	\begin{figure}[!t]
%	\centering
%	\includegraphics[width=2.5in]{fig_f}
%	\caption{Plotted here is the roll-off slope of the hyperbolic window at $\alpha=3$, contrasted with that of the fixed-width windows (open circles), the Tukey (dotted line) and the Planck (dashed line). The solid line gives the roll-off calculated from the total curve; the dot-dashed line shows the roll-off calculated from the residual curve. The difference between these is described in section \ref{ss_hyperbolic}. Note the discontinuity in the `total curve' line at the Hann-cubed value. As explained in section \ref{ss_hyperbolic}, this is due to the rapid growth of a non-oscillatory perturbation to the FT as the warp factor moves away from zero.}
%	\label{fig_f}
%	\end{figure}

The slow decrease in the steepness of the roll-off for the hyperbolic window seen in Fig. \ref{fig_f} at ENBW less than about 1.2 is also an artifact of the empirical $Tf \sim 150$ estimation point. The hyperbolic window in FT exhibits a pedestal which becomes ever broader as ENBW is decreased. The estimates of roll-off presented here thus begin to decrease as the width of the pedestal approaches the $Tf = 150$ value.

	\begin{figure}[!t]
	\centering
	\includegraphics[width=2.5in]{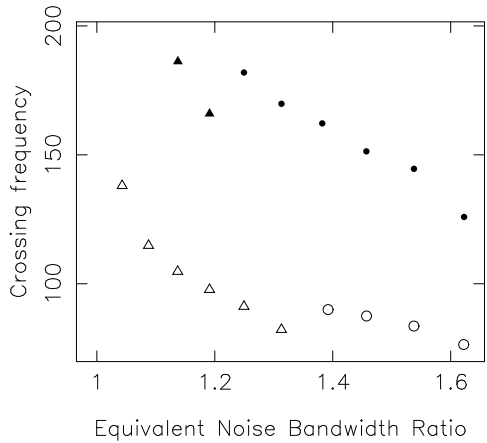}
	\caption{This Fig. shows the approximate frequency at which the Fourier transform of the Planck window of the given ENBW first becomes smaller than the transform of the hyperbolic window of the same ENBW. As with Fig. \ref{fig_i}, triangles and circles indicate respectively where the total or residual curves of the hyperbolic FT were considered. Open symbols indicate $\alpha=3$; filled ones, $\alpha=4$.}
	\label{fig_h}
	\end{figure}

In fact the Planck window FT also exhibits a pedestal which is generally seen to be wider than that of the hyperbolic window at the same ENBW. The respective FTs don't cross until typically $Tf \sim 100$. Fig. \ref{fig_h} shows the (approximate) frequency at which the Planck FT first dips below the hyperbolic FT, for two different values of $\alpha$. For all frequencies below this value, the hyperbolic FT is of lower amplitude than the Planck, by as much as 30 decibels. Broadly speaking we find that, for a given value of ENBW, an exponent $\alpha$ can be found at which the hyperbolic is a better performer than the Planck up to frequencies of about $100/T$; beyond that point, the Planck is usually better. With the double-precision arithmetic on the machine used to perform the calculations, the noise floor of about $-340$ decibels is reached by the hyperbolic window FT at $Tf$ values of order 1000.

%++++++++++++++++++++++++++++++++++++++++++++++++++++++++++++
\section{Conclusions}

The present paper recaps the basic considerations behind the use of window functions to assist discrete Fourier analysis of signals comprising one or more periodic components among noise. Of the many possible ways to assess the performance of windows, attention is concentrated here on those criteria which bear directly on either sensitivity or dynamic range. Three quantitative criteria were selected, namely the equivalent noise bandwidth or ENBW, the amplitude of the highest sidelobe of the window Fourier transform, and the high-frequency asymptotic rate of decrease of the transform (roll-off).

Although a low value of ENBW is associated not only with optimum sensitivity but also with a narrower central maximum in the window FT, there is no discussion of resolution here. Coverage of all the different criteria is beyond the scope of the present paper. Maximizing the resolving power is in any case more complicated than simply choosing the narrowest central peak: high sidelobes will play a role as well.

A new variable-width window, named the hyperbolic window, is described in this paper. This window, at representative values of its width and power parameters $s$ and $\alpha$, is compared against a selection of windows from the existing literature. Since, by choice of the appropriate value of $s$, the ENBW of the hyperbolic window can be varied continuously between its minimum and maximum possible values of respectively unity and infinity, the most important comparisons display either sidelobe height or roll-off for the panel of windows as smooth functions of ENBW. The effect of $\alpha$ on the performance of the hyperbolic window is indicated where practical by plotting curves at a few values of $\alpha$ on the same graph. Note that, since best sensitivity is obtained for low values of ENBW, attention has been concentrated on this range.

Unlike the Tukey or Planck windows, the width of the hyperbolic window may be decreased without limit: in other words it can access all values of ENBW up to infinity. However since there seems nothing to be gained by going to large ENBW, this regime is little explored in the present paper.

Because the Planck window is smooth at all orders of derivative, the steepness of its FT roll-off increases with frequency without bound. In practice however, although this gives it superior dynamic range over a range of frequencies, its wide pedestal sets a lower boundary to this range at $Tf \sim 100$, and the finite numerical precision of computation sets an upper one at $Tf \sim 1000$.

Another advantage the hyperbolic window has over both the Tukey and Planck windows is that it is not constructed by splicing two different functions together. This gives it a much more regular-appearing Fourier transform. Both the Tukey and Planck windows in contrast exhibit relatively cluttered spectra, which may marginally decrease the detectability of weak additional tones.

A feature of the hyperbolic window is that, for significant ranges of $s$ and $\alpha$, its Fourier transform ceases at some frequency to be oscillatory, i.e. to pass through zero. (Indeed for some values of $s$ and $\alpha$ there are no zero-crossings at all.) In such cases it is convenient to divide the transform into a DC offset curve plus an oscillatory residual. Usually the oscillatory part has a much steeper roll-off than the smooth part. It is shown that, as far as detecting secondary tones goes, it is the oscillatory not the smooth part which is the practical limiting factor.

The comparison with other windows shows that, for any ENBW and at any frequency less than about $100/T$, the hyperbolic window is better either in sidelobe height or roll-off than any of its competitors - although not necessarily in both. It never approaches the sidelobe performance of the Nuttall windows at the same ENBW values, although for $\alpha > 2$ it has better roll-off. Nevertheless, it offers a respectable 40 decibels of sidelobe suppression at an ENBW of about 1.5, whereas for the top-performing Nuttall windows ENBW is greater than 2.

\bibliographystyle{plainnat}
\bibliography{./warp}

\end{document}